\documentclass[11pt]{article}
\usepackage[dvipsnames]{xcolor}
\usepackage{amsfonts}
\usepackage{setspace}
\usepackage{bm}
\usepackage{bbm}
\usepackage{amssymb}
\usepackage{indentfirst}
\usepackage{geometry}
\usepackage{mathtools}
\usepackage{hyperref}
\hypersetup{colorlinks=false, hidelinks, urlcolor=cyan}

\DeclarePairedDelimiter\abs{\lvert}{\rvert}%
\DeclarePairedDelimiter{\floor}{\lfloor}{\rfloor}
\DeclareMathOperator{\HP}{HP}

\DeclareMathOperator{\csch}{csch}

\makeatletter
\let\oldabs\abs
\def\abs{\@ifstar{\oldabs}{\oldabs*}}

\newcommand*{\citena}[1]{%
\begingroup
\color{Green}
\romannumeral-`\x 
\setcitestyle{numbers}%
\citep{#1}
\endgroup
\ignorespacesafterend
}

\newcommand*{\citesup}[1]{%
\begingroup
\color{Green}
\citep{#1}
\endgroup
\ignorespacesafterend
}

\newcommand*{\eqrefe}[1]{%
\begingroup
(\color{BrickRed}
\romannumeral-`\x 
\setcitestyle{numbers}%
\ref{eq:#1}%
\endgroup
)\ignorespacesafterend
}

\newcommand*{\secrefe}[1]{%
\begingroup
(\color{Aquamarine}
\romannumeral-`\x 
\setcitestyle{numbers}%
\ref{#1}%
\endgroup
)\ignorespacesafterend
}

\newtagform{Tags}[\textcolor{BrickRed}]{\color{Black}(}{)}

\usepackage[numbers, sort]{natbib}
\setcitestyle{super}

\newcommand{\ii}{\bm{i}}

\begin{document}
\title{Generalized Harmonic Progression Part II}
\author{
Jose Risomar Sousa\\
\href{https://orcid.org/0009-0002-7779-0808}
     {\textcolor{cyan}{ORCID: 0009-0002-7779-0808}}
}
\date{February 4, 2019}
\maketitle

\begin{abstract}
In a previous paper, we saw how to create formulae for the sum of the terms of a harmonic progression of order $k$, $\HP_k(n)$, with integer parameters, $a$ and $b$. In this new paper we make those formulae even more general by lifting the restriction that the parameters be integers. The new formula holds always, except when $\ii\,b/a \in \mathbb{Z}$. Here a slightly modified version of the reasoning used before is introduced, in a standalone exposition that does not require prior knowledge of the precursor paper. A notable advantage of this new approach is that it can be used to find summation formulae for the reciprocal of polynomials, that is, $\sum_{j}1/p(j)$.
\end{abstract}

\usetagform{Tags}

\tableofcontents


\section{Introduction}

This manuscript demonstrates how to produce formulae for a harmonic progression of order $k$ with complex parameters,
\begin{equation} \nonumber
\HP_k(n)=\sum_{j=1}^{n}\frac{1}{(a\ii j+b)^k} \text{,}
\end{equation}
\noindent with $a$ an integer, $b$ a complex number and $\ii$ the imaginary unit.\\

The motivation to create such a formula came from the realization that, unlike the sum of $1/(j^2-1)$ over $j$, a sum such as $1/(j^2+1)$ can not be derived from the formulae created in \citena{GHP}.\\ 

If a formula for $\HP_k(n)$ with integer $a$ and complex $b$ can be produced, formulae for most, if not all, sums of the type $\sum_{j}1/p(j)$, with $p(j)$ a polynomial in $j$ with complex coefficients, can be produced. Besides, $\HP_k(n)$ provides the partial sums of the Hurwitz zeta function, $\zeta(k,b)$, at the positive integers $k$, whose formula can be obtained by taking the limit of $\HP_k(n)$ when $n$ goes to infinity. With a little more effort, it may even be possible to obtain the analytic continuation of $\zeta(k,b)$.\\

The choice of integer $a$ and complex $b$ is good enough in most of the cases (since one can always make $a$ integer by dividing the parameters by $a$), and when it is not, one can resort to the formulae from the previous paper\citesup{GHP}, which assume $a$ and $b$ to be integers.\\

The method makes use again of the Faulhaber formula\citesup{John} for the sum of the powers of the first $n$ positive integers, which are given by,
\begin{equation} \label{eq:soma_pot_par}
\sum_{j=1}^{n} {j^{2i}}=\frac{n^{2i}}{2}+\sum_{j=0}^{i} \frac{(2i)! B_{2j}\,n^{2i+1-2j}}{(2j)!(2i+1-2j)!} 
\end{equation} 
\begin{equation} \label{eq:soma_pot_impar}
\sum_{j=1}^{n} {j^{2i+1}}=\frac{n^{2i+1}}{2}+\sum_{j=0}^{i} \frac{(2i+1)! B_{2j}\, n^{2i+2-2j}}{(2j)!(2i+2-2j)!} \text{,}
\end{equation}
\noindent where $B_{2j}$ are the non-zero Bernoulli numbers.

\section{Review of previous results}

Before we begin, let us recapitulate some results from \citena{GHP} for a bit of context. We had seen that the sum of the terms of a generalized harmonic progression with integer parameters, $a$ and $b$, for even and odd powers is given, respectively, by,
\begin{multline} \nonumber 
\sum_{j=1}^{n}\frac{1}{(a j+b)^{2k}}=-\frac{1}{2b^{2k}}+\frac{1}{2(a n+b)^{2k}}
\\-\frac{(-1)^{k}(2\pi)^{2k}}{2}\int_{0}^{1}\sum_{j=0}^{k}\frac{B_{2j}\left(2-2^{2j}\right)(1-u)^{2k-2j}}{(2j)!(2k-2j)!}\left(\sin{2\pi(a n+b)u}-\sin{2\pi bu}\right)\cot{\pi a u}\,du \text{,}
\end{multline}
\noindent and
\begin{multline} \nonumber 
\sum_{j=1}^{n}\frac{1}{(a j+b)^{2k+1}}=-\frac{1}{2b^{2k+1}}+\frac{1}{2(a n+b)^{2k+1}}\\-\frac{(-1)^{k}(2\pi)^{2k+1}}{2}\int_{0}^{1}\sum_{j=0}^{k}\frac{B_{2j}\left(2-2^{2j}\right)(1-u)^{2k+1-2j}}{(2j)!(2k+1-2j)!}\left(\cos{2\pi(a n+b)u}-\cos{2\pi bu}\right)\cot{\pi au}\,du \text{}
\end{multline}
\indent Notice that there may be singularities on both sides of these equations if, for example, $a j+b=0$ for some $j\leq n$ (that is, $a$ divides $b$ and they have different signs). In those cases, if the singularities are ignored wherever they occur the formula still holds.\\ 

This means that, for instance, if $b=0$ then the above formulae reduce to the formulae of the generalized harmonic numbers from reference \citena{GHN}.\\

To see why singularities can be removed, let us look at the formula for the odd powers. From that formula, it is easy to see why it still holds when singularities are ignored. For example, let us take the forward difference of $\HP(n)$,
\begin{multline} \nonumber
2\pi\int_{0}^{1}(1-u)\left(\cos{2\pi(a n+b)u}-\cos{2\pi(a(n-1)+b)u}\right)\cot{\pi a(1-u)}\,du \\ =-\frac{1}{a n+b}-\frac{1}{a(n-1)+b}\text{}
\end{multline}
\indent Now, if for instance $a n+b=0$, the left-hand side of the above equation yields $1/a$, even though $1/(a n+b)$ is a singularity on the right-hand side. Therefore, by ignoring the singularity, the equation still holds.\\

There was also a unified formula based on the exponential (that is, the initial equation). For every integer $a$, $b$ and integer $k\ge 1$:
\begin{multline} \label{eq:expo_int_pmts}
\sum_{j=1}^{n}\frac{1}{(a j+b)^{k}}=-\frac{1}{2b^{k}}+\frac{1}{2(a n+b)^{k}}\\+\frac{\ii(2\pi\ii)^{k}}{2}\int_{0}^{1}\sum_{j=0}^{k}\frac{B_{j}(1-u)^{k-j}}{j!(k-j)!}\left(e^{2\pi\ii(a n+b)u}-e^{2\pi\ii b\,u}\right)\cot{\pi a u}\,du 
\end{multline}

\section{Approach based on exponential} \label{HN_sin_k_pi}

The reasoning to build a formula for the harmonic progression relies on the Taylor series expansion of the exponential, and seizes upon the fact that it is constant when $a$ is integer and $b$ is complex:
\begin{equation} \label{eq:start_iden} 
e^{2\pi (a\ii j+b)}=e^{2\pi b} \Rightarrow  \sum_{i=0}^{\infty}\frac{(2\pi (a\ii j+b))^i}{i!}=e^{2\pi b} \text{, if }a\text{ is integer} 
\end{equation}
\indent Starting from the above initial equation, one ends up with the below after a few steps:
\begin{equation} \label{eq:start_iden_x}
\frac{1}{a\ii j+b}(e^{2\pi b}-1)=2\pi\sum_{i=0}^{\infty}\left(\frac{(2\pi(a\ii j+b))^{2i}}{(2i+1)!}+\frac{(2\pi(a\ii j+b))^{2i+1}}{(2i+2)!}\right)
\end{equation}
\indent Note there are now two power series, as opposed to one previously (that is, in \citena{GHN} and \citena{GHP}).\\

At this point it is evident that this formula does not work if $e^{2\pi b}=1$, which happens when $b$ is a purely complex number with integer imaginary part ($\Re(b)=0$ and $\Im(b)\in \mathbb{Z}$). In most situations, this issue can be avoided by making $b$ non-integer, such as by dividing $b$ by $a$. Likewise, the formula does not apply if $a$ is not integer ($a\notin \mathbb{Z}$). In that case $a$ can be made integer by dividing both parameters by $a$:
\begin{equation} \nonumber
\sum_{j=1}^{n}\frac{1}{a\ii j+b}=\frac{1}{a}\sum_{j=1}^{n}\frac{1}{\ii j+b/a} 
\end{equation}
\indent However, in both cases this solution does not work if $\ii\,b/a$ is integer, including zero. In those cases the previous formulae are applicable.

\subsection{Lagrange's trigonometric identities} \label{Lagra}

Through Lagrange's identities, which stem from the sum of the terms of a geometric progression, functions can be created whose power series take similar forms to the power series that are obtained when the initial equation from the previous section is expanded.\\

For the complex harmonic progression, these functions are the same as the ones used in \citena{GHP}, with only a slight transformation ($a \mapsto a\ii$). Nonetheless, they are obtained by means of the Lagrange's trigonometric identities, 
\begin{multline} \label{eq:cos} 
\sum_{j=1}^{k}\cos{\frac{2\pi n(a\ii j+b)}{k}}=-\frac{1}{2}\cos{\frac{2\pi b n}{k}}+\frac{1}{2}\cos{2\pi n\left(a\ii+\frac{b}{k}\right)}\\+\cos{\pi n\left(a\ii+\frac{2b}{k}\right)}\sin{\pi a\ii n}\cot{\frac{\pi a\ii n}{k}} \text{,}
\end{multline}
\noindent and
\begin{multline} \label{eq:seno} 
\sum_{j=1}^{k}\sin{\frac{2\pi n(a\ii j+b)}{k}}=-\frac{1}{2}\sin{\frac{2\pi b n}{k}}+\frac{1}{2}\sin{2\pi n\left(a\ii+\frac{b}{k}\right)}\\+\sin{\pi n\left(a\ii+\frac{2b}{k}\right)}\sin{\pi a\ii n}\cot{\frac{\pi a\ii n}{k}} \text{,}
\end{multline}
\noindent along with the cosine and sine addition formulae\citesup{Abra}:
\begin{equation} \nonumber
\cos(x+y)=\cos{x}\cos{y}-\sin{x}\sin{y} \text{, and } \sin(x+y)=\sin{x}\cos{y}+\cos{x}\sin{y} 
\end{equation}
\indent Power series for the left-hand side of the Lagrange identities can be derived with the employment of \eqrefe{soma_pot_par} and \eqrefe{soma_pot_impar}, resulting in the following power series for each function on the right-hand side:
\begin{equation} \nonumber
\cos{2\pi n\left(a\ii+\frac{b}{k}\right)}=\sum_{i=0}^{\infty}(-1)^i\left(\frac{2\pi b n}{k}\right)^{2i}\left(\sum_{j=0}^{i}\frac{(a\ii k/b)^{2j}}{(2j)!(2i-2j)!}+\sum_{j=0}^{i-1}\frac{(a\ii k/b)^{2j+1}}{(2j+1)!(2i-2j-1)!}\right)
\end{equation} 

\begin{equation} \nonumber
\sin{2\pi n\left(a\ii+\frac{b}{k}\right)}=\sum_{i=0}^{\infty}(-1)^i\left(\frac{2\pi b n}{k}\right)^{2i+1}\left(\sum_{j=0}^{i}\frac{(a\ii k/b)^{2j}}{(2j)!(2i+1-2j)!}+\sum_{j=0}^{i}\frac{(a\ii k/b)^{2j+1}}{(2j+1)!(2i-2j)!}\right)
\end{equation}

\begin{multline} \nonumber
\cos{\pi n\left(a\ii+\frac{2b}{k}\right)}\sin{\pi a\ii n}\cot{\frac{\pi a\ii n}{k}}=\\\sum_{i=0}^{\infty}(-1)^i\left(\frac{2\pi b n}{k}\right)^{2i}\left(\sum_{j=0}^{i}\frac{(a\ii/b)^{2j}}{(2i-2j)!}\sum_{q=0}^{j}\frac{B_{2q}k^{2j+1-2q}}{(2j+1-2q)!(2q)!}+\sum_{j=0}^{i-1}\frac{(a\ii /b)^{2j+1}}{(2i-2j-1)!}\sum_{q=0}^{j}\frac{B_{2q}k^{2j+2-2q}}{(2j+2-2q)!(2q)!}\right)
\end{multline}

\begin{multline} \nonumber
\sin{\pi n\left(a\ii+\frac{2b}{k}\right)}\sin{\pi a\ii n}\cot{\frac{\pi a\ii n}{k}}=\\\sum_{i=0}^{\infty}(-1)^i\left(\frac{2\pi b n}{k}\right)^{2i+1}\left(\sum_{j=0}^{i}\frac{(a\ii/b)^{2j}}{(2i+1-2j)!}\sum_{q=0}^{j}\frac{B_{2q}k^{2j+1-2q}}{(2j+1-2q)!(2q)!}+\sum_{j=0}^{i}\frac{(a\ii /b)^{2j+1}}{(2i-2j)!}\sum_{q=0}^{j}\frac{B_{2q}k^{2j+2-2q}}{(2j+2-2q)!(2q)!}\right)
\end{multline}

\subsection{Harmonic progression} \label{HN_1_sin_2k_pi}

To obtain $\HP(n)$, equation \eqrefe{start_iden_x} is summed over $j$, and then $(a\ii j+b)^{2i}$ and $(a\ii j+b)^{2i+1}$ are expanded using the binomial theorem, and lastly the sums of $j^{2q}$ and $j^{2q+1}$ over $j$ are replaced with their respective Faulhaber formulae. After all the possible simplifications are carried out, one ends up with a few power series, shown below.\\

First, one has the independent term given by:
\begin{equation} \nonumber
-\pi\sum_{i=0}^{\infty}\left(\frac{(2\pi b)^{2i}}{(2i+1)!}+\frac{(2\pi b)^{2i+1}}{(2i+2)!}\right)=-\frac{e^{2\pi b}-1}{2b}
\end{equation}
\indent For the other, more complicated power series, integral representations can be obtained with the aid of their corresponding functions from section \secrefe{Lagra}:
\begin{equation} \nonumber
\pi\sum_{i=0}^{\infty}\frac{(2\pi b)^{2i}}{2i+1}\left(\sum_{j=0}^{i}\frac{(a\ii n/b)^{2j}}{(2j)!(2i-2j)!}+\sum_{j=0}^{i-1}\frac{(a\ii n/b)^{2j+1}}{(2j+1)!(2i-2j-1)!}\right)=-\frac{\ii \pi}{n}\int_{0}^{\ii n}\cos{2\pi x\left(a\ii+\frac{b}{n}\right)}\,dx 
\end{equation} 

\begin{equation} \nonumber
\pi\sum_{i=0}^{\infty}\frac{(2\pi b)^{2i+1}}{2i+2}\left(\sum_{j=0}^{i}\frac{(a\ii n/b)^{2j}}{(2j)!(2i+1-2j)!}+\sum_{j=0}^{i}\frac{(a\ii n/b)^{2j+1}}{(2j+1)!(2i-2j)!}\right)=-\frac{\pi}{n}\int_{0}^{\ii n}\sin{2\pi x\left(a\ii+\frac{b}{n}\right)}\,dx 
\end{equation} 

\begin{multline} \nonumber
2\pi\sum_{i=0}^{\infty}\frac{(2\pi b)^{2i}}{2i+1}\left(\sum_{j=0}^{i}\frac{(a\ii/b)^{2j}}{(2i-2j)!}\sum_{q=0}^{j}\frac{B_{2q}n^{2j+1-2q}}{(2j+1-2q)!(2q)!}+\sum_{j=0}^{i-1}\frac{(a\ii/b)^{2j+1}}{(2i-2j-1)!}\sum_{q=0}^{j}\frac{B_{2q}n^{2j+2-2q}}{(2j+2-2q)!(2q)!}\right) 
\\=-\frac{2\pi\ii}{n}\int_{0}^{\ii n}\cos{\pi x\left(a\ii+\frac{2b}{n}\right)}\sin{(\pi a\ii x)}\cot{\frac{\pi a\ii x}{n}}\,dx 
\end{multline}

\begin{multline} \nonumber
2\pi\sum_{i=0}^{\infty}\frac{(2\pi b)^{2i+1}}{2i+2}\left(\sum_{j=0}^{i}\frac{(a\ii/b)^{2j}}{(2i+1-2j)!}\sum_{q=0}^{j}\frac{B_{2q}n^{2j+1-2q}}{(2j+1-2q)!(2q)!}+\sum_{j=0}^{i}\frac{(a\ii/b)^{2j+1}}{(2i-2j)!}\sum_{q=0}^{j}\frac{B_{2q}n^{2j+2-2q}}{(2j+2-2q)!(2q)!}\right)\\=-\frac{2\pi}{n}\int_{0}^{\ii n}\sin{\pi x\left(a\ii+\frac{2b}{n}\right)}\sin{(\pi a\ii x)}\cot{\frac{\pi a\ii x}{n}}\,dx 
\end{multline}
\indent Therefore, after putting it all together and making a change of variables, the formula below is found, which holds as long as $a$ is integer and $\ii\,b$ is not integer:
\begin{equation} \label{eq:HP(n)_1}
\sum_{j=1}^{n}\frac{1}{a\ii j+b}=-\frac{1}{2b}+\frac{1}{2(a\ii n+b)}+\frac{2\pi}{e^{2\pi b}-1}\int_{0}^{1}e^{\pi(a\ii n+2b)u}\sin{\pi a n u}\cot{\pi a u}\,du 
\end{equation}
\indent In this formula, the independent term came from the first power series, the second term from the next two power series and the third term from the last two power series, which is somewhat evident.\\

As discussed before, provided that $\ii\,b/a\notin \mathbb{Z}$, the restriction that $a$ be integer can be lifted by transforming the formula, for example, as follows,
\begin{equation} \label{eq:HP(n)_2}
\frac{1}{a}\sum_{j=1}^{n}\frac{1}{\ii j+b/a} =-\frac{1}{2b}+\frac{1}{2(a\ii n+b)}+\frac{2\pi}{a(e^{2\pi b/a}-1)}\int_{0}^{1}e^{\pi(\ii n+2b/a)u}\sin{\pi n u}\cot{\pi u}\,du \text{,} 
\end{equation}
\noindent which holds for any $a,$ $b\in\mathbb{C}$. There are endless ways of doing this, this example is just one of the possibilities.

\subsection{Generalization} \label{General}

To generalize the previous formula to $\HP_k(n)$, first note that a recurrence is obtained at each new step. For example, for the third power one has:
\begin{multline} \nonumber 
(e^{2\pi b}-1)\sum_{j=1}^{n}\frac{1}{(a\ii j+b)^3}=\\\sum_{j=1}^{n}\left(\frac{(2\pi)^2}{2!}\frac{1}{a\ii j+b}+\frac{2\pi}{1!}\frac{1}{(a\ii j+b)^2}+(2\pi)^3\sum_{i=0}^{\infty}\frac{(2\pi(a\ii j+b))^{2i}}{(2i+3)!}+\frac{(2\pi(a\ii j+b))^{2i+1}}{(2i+4)!}\right)
\end{multline}
\indent This means that only the last sum in the equation above is complicated. As before, this recurrence is such that in the final formula the terms that go outside of the integral reduce to very simple forms, with the integral being the only challenging part.\\

\indent Let $p_k(u)$ be the polynomial in $u$ that goes within the integral. Then for integer $a$ and non-integer $\ii\,b/a$,
\begin{equation} \label{eq:general_form}
\sum_{j=1}^{n}\frac{1}{(a\ii j+b)^k}=-\frac{1}{2b^k}+\frac{1}{2(a\ii n+b)^k}+(2\pi)^k\int_{0}^{1}p_k(u)e^{\pi(a\ii n+2b)u}\sin{\pi a n u}\cot{\pi a u}\,du \text{,}
\end{equation}
\noindent where  $p_k(u)$ is given by the recurrence equation,
\begin{equation} \nonumber
(e^{2\pi b}-1)p_{k}(u)=
\begin{cases}
      1, & \text{if}\ k=1\\
      \frac{(1-u)^{k-1}}{(k-1)!}+\sum_{j=1}^{k-1}\frac{p_j(u)}{(k-j)!}, & \text{if integer }\ k>1
\end{cases} 
\end{equation}
and where term $(1-u)^{k-1}/(k-1)!$ came from the observation of the patterns.\\

Now to solve the recurrence, let $p(x)$ be the generating function of $p_k(u)$. By analyzing the recurrence, it is possible to conclude that:
\begin{equation} \label{eq:recur}
(e^{2\pi b}-1)p(x)-\left(-1+e^x\right)p(x)=x\,e^{(1-u)x}\Rightarrow p(x)=-\frac{x\,e^{(1-u)x}}{e^{x}-e^{2\pi b}}
\end{equation}
\indent The solution of this equation is not trivial, as the general term of the power series of $p(x)$ is not simple.\\

To find $p_k(u)$, let $\delta_{ij}$ be the Kronecker delta ($\delta_{ij}=1$, if $i=j$, $\delta_{ij}=0$ otherwise), and $\mathrm{Li}_{k}(z)$, the polylogarithm\citesup{Abra}.\\

The $p_k(u)$ are the coefficients of the Taylor series expansion of $p(x)$, which are,
\begin{equation} \nonumber
p_k(u)=\frac{p^{(k)}(0)}{k!}=\sum_{j=1}^{k}\frac{\left(\delta_{1j}+\mathrm{Li}_{-j+1}(e^{-2\pi b})\right)(1-u)^{k-j}}{(j-1)!(k-j)!} \text{,}
\end{equation}
\noindent and the identity below,
\begin{equation} \nonumber
\sum_{j=1}^{k}\frac{\mathrm{Li}_{-j+1}(e^{-2\pi b})(1-u)^{k-j}}{(j-1)!(k-j)!}=\frac{e^{-2\pi b}}{\Gamma(k)}\,\Phi(e^{-2\pi b},-k+1,2-u) \text{,}
\end{equation}
\noindent provides an extension of the summation to non-integer $k$.\\

Further simplifications can be attained through the identity below,
\begin{equation} \label{eq:Lerch_neg_id}
(1-u)^{k-1}+e^{-2\pi b}\,\Phi(e^{-2\pi b},-k+1,2-u)=\Phi(e^{-2\pi b},-k+1,1-u) \text{,}
\end{equation}
\noindent that stems from the definition of the Lerch $\Phi$ function and also through Euler's formula:
\begin{equation} \label{eq:Euler_form}
e^{\pi\ii a n u}\sin{\pi a n u}=(\cos{\pi a n u}+\ii\,\sin{\pi a n u})\sin{\pi a n u}=\frac{\ii}{2}\left(1-e^{2\pi\ii a n u}\right)\text{,}
\end{equation}
\indent Therefore, if $k$ is complex with $\Re(k)>0$ and $a$ is integer, one has,
\begin{multline} \label{eq:Extended_form_1_orig}
\sum_{j=1}^{n}\frac{1}{(a\ii j+b)^k} =-\frac{1}{2b^k}+\frac{1}{2(a\ii n+b)^k}\\-\frac{\ii\,(2\pi)^k}{2\,\Gamma(k)}\int_{0}^{1}\,e^{-2\pi b(1-u)}\,\Phi(e^{-2\pi b},-k+1,1-u)\,(e^{2\pi\ii a n u}-1)\cot{\pi a u}\,du 
\end{multline}
\indent Finally, the formula can be turned into a more straightforward form by making the offset parameter equal to $\ii\,b$:
\begin{multline} \label{eq:Extended_form_2_orig}
\sum_{j=1}^{n}\frac{1}{(a j+b)^k} =-\frac{1}{2b^k}+\frac{1}{2(a n+b)^k}\\-\frac{\ii\,(2\pi\ii)^k}{2\,\Gamma(k)}\int_{0}^{1}\,e^{-2\pi\ii\,b(1-u)}\,\Phi(e^{-2\pi\ii\,b},-k+1,1-u)\,(e^{2\pi\ii a n u}-1)\cot{\pi a u}\,du \text{}
\end{multline}

\subsection{Solution of the transform}

A transformation of \eqrefe{general_form} is shown below, which is just setting $a=1$ and $b \mapsto b/a$ and dividing both sides by $a^{k}$, which can be more useful for the purpose of finding the limit of the harmonic progression of order $k$:
\begin{equation} \label{eq:general_transform}
\sum_{j=1}^{n}\frac{1}{(a\ii j+b)^k}=-\frac{1}{2b^k}+\frac{1}{2(a\ii n+b)^k}+\left(\frac{2\pi}{a}\right)^k\int_{0}^{1}p_k(u)e^{\pi(\ii n+2b/a)u}\sin{\pi n u}\cot{\pi u}\,du \text{}
\end{equation}
\indent Now, the generating function of the $p_k(u)$ becomes:
\begin{equation} \nonumber
p(x)=-\frac{x\,e^{x(1-u)/a}}{a(e^{x/a}-e^{2\pi b/a})} \Rightarrow \frac{p^{(k)}(0)}{k!}=\frac{e^{-2\pi b/a}}{a^{k}(k-1)!}\,\Phi(e^{-2\pi b/a},-k+1,1-u) \text{,}
\end{equation}\\
\indent Therefore, if $\Re(k)>0$ and $\ii\,b/a\notin \mathbb{Z}$, one has:
\begin{multline} \label{eq:Extended_form_1}
\sum_{j=1}^{n}\frac{1}{(a\ii j+b)^k} =-\frac{1}{2b^k}+\frac{1}{2(a\ii n+b)^k}\\-\frac{\ii}{2\,\Gamma(k)}\left(\frac{2\pi}{a}\right)^k\int_{0}^{1}\,e^{-2\pi(1-u)b/a}\,\Phi(e^{-2\pi b/a},-k+1,1-u)\,(e^{2\pi\ii n u}-1)\cot{\pi u}\,du 
\end{multline}
\indent The formula can be simplified further by letting $a=-\ii$. Thus, if $\Re(k)>0$ and $b\notin \mathbb{Z}$ then:
\begin{multline} \label{eq:Extended_form_2}
\sum_{j=1}^{n}\frac{1}{(j+b)^k} =-\frac{1}{2b^k}+\frac{1}{2(n+b)^k}\\-\frac{\ii\,(2\pi\ii)^k}{2\,\Gamma(k)}\int_{0}^{1}\,e^{-2\pi\ii b(1-u)}\,\Phi(e^{-2\pi\ii\,b},-k+1,1-u)\,(e^{2\pi\ii n u}-1)\cot{\pi u}\,du \text{}
\end{multline}

\subsection{Generalization to reciprocal of polynomials} \label{Examples}

Obtaining formulae for the sum of the reciprocal of a polynomial $p(j)$ of order $k$ over $j$ is straightforward, as long as the decomposition of the polynomial does not have any multiple roots. If $x_q$, $q=1, \ldots, k$, are the $k$ unique roots of $p(j)$ then:
\begin{equation}
\frac{1}{p(j)}=\sum_{q=1}^{k}\frac{c_{q}}{j-x_q}
\end{equation}
\indent If any of the roots is multiple, then the decomposition is a little more complicated, but the literature\citesup{Apostol} has a general method to find the $c_q$ as a function of $p(j)$, its roots $x_q$ and their multiplicity, anyway.\\

As a particular case, to obtain a formula for the sum of $1/(j^{2k}+1)$ over $j$ it is necessary to find the $2k$ complex roots of $x^{2k}+1=0$, as below:
\begin{equation} \label{eq:form_gen}
\sum_{q=1}^{2k}\frac{c_{q}}{j-x_q}=\frac{1}{j^{2k}+1} \text{}
\end{equation}
\indent Let us see two examples on how the newly found formulae can be applied to obtain sums that were not possible with the previous formulae created in \citena{GHP}.\\

The first example is the sum of $1/(j^{2}+1)$ over $j$. The first step to derive this formula is to decompose it into a linear combination as below:
\begin{equation} \nonumber
\frac{1}{j^2+1}=\frac{1}{2(\ii j+1)}-\frac{1}{2(\ii j-1)}
\end{equation}
\indent After a few simplifications, the formula below is found:
\begin{equation} \label{eq:form_1}
\sum_{j=1}^{n}\frac{1}{j^2+1}=-\frac{1}{2}+\frac{1}{2(n^2+1)}+\pi\csch{2\pi}\int_{0}^{1}\cosh{2\pi(2u-1)}\sin{4\pi n u}\cot{2\pi u}\,du
\end{equation}

\indent The second example is another quadratic polynomial, which can be decomposed as below,
\begin{equation} \nonumber
\frac{1}{j^2+2j+2}=\frac{1}{2(\ii j+1+\ii)}-\frac{1}{2(\ii j-1+\ii)}\text{,}
\end{equation}
\noindent and which therefore implies the formula:
\begin{multline} \label{eq:form_2}
\sum_{j=1}^{n}\frac{1}{j^2+2j+2}=-\frac{1}{4}+\frac{1}{2(n^2+2n+2)}\\+\frac{2\pi}{e^{4\pi}-1}\int_{0}^{1}(e^{4\pi(1-u)}+e^{4\pi u})\cos{2\pi(n+2)u}\sin{2\pi n(1-u)}\cot{2\pi(1-u)}\,du
\end{multline}
\indent Therefore, it becomes evident that the real problem is finding the linear combination that yields the reciprocal of the polynomial that one is interested in. 

\section{Approaches based on cosine and sine}
Though the exponential-based approach is probably more useful, if $\cos{2\pi(j+b)}=\cos{2\pi b}$ and $\sin{2\pi(j+b)}=\sin{2\pi b}$ are used as initial identities, different formulae can be obtained.\\

For the cosine, starting from,
\begin{equation} \nonumber
\cos{2\pi(j+b)}=\cos{2\pi b}=\sum_{i=0}^{\infty}\frac{(-1)^i(2\pi(j+b))^{2i}}{(2i)!} \text{,}
\end{equation}
\noindent produces the recurrence equations:
\begin{equation} \nonumber
\begin{cases}
\frac{2(\sin{\pi b})^2}{j+b}=2\pi\sum_{i=0}^{\infty}\frac{(-1)^i(2\pi(j+b))^{2i+1}}{(2i+2)!}\\

\frac{2(\sin{\pi b})^2}{(j+b)^2}=(2\pi)^2\sum_{i=0}^{\infty}\frac{(-1)^i(2\pi(j+b))^{2i}}{(2i+2)!}\\

\frac{2(\sin{\pi b})^2}{(j+b)^3}=\frac{(2\pi)^2}{2!(j+b)}-(2\pi)^3\sum_{i=0}^{\infty}\frac{(-1)^i(2\pi(j+b))^{2i+1}}{(2i+4)!}\\

\frac{2(\sin{\pi b})^2}{(j+b)^4}=\frac{(2\pi)^2}{2!(j+b)^2}-(2\pi)^4\sum_{i=0}^{\infty}\frac{(-1)^i(2\pi(j+b))^{2i}}{(2i+4)!}\\

\frac{2(\sin{\pi b})^2}{(j+b)^5}=\frac{(2\pi)^2}{2!(j+b)^3}-\frac{(2\pi)^4}{4!(j+b)}+(2\pi)^5\sum_{i=0}^{\infty}\frac{(-1)^i(2\pi(j+b))^{2i+1}}{(2i+6)!}\\

\frac{2(\sin{\pi b})^2}{(j+b)^6}=\frac{(2\pi)^2}{2!(j+b)^4}-\frac{(2\pi)^4}{4!(j+b)^2}+(2\pi)^6\sum_{i=0}^{\infty}\frac{(-1)^i(2\pi(j+b))^{2i}}{(2i+6)!}\\

\vdots 

\end{cases}
\end{equation}
\noindent Likewise, for the sine, starting from,
\begin{equation} \nonumber
\sin{2\pi(j+b)}=\sin{2\pi b}=\sum_{i=0}^{\infty}\frac{(-1)^i(2\pi(j+b))^{2i+1}}{(2i+1)!} \text{,}
\end{equation}
\noindent the recurrence equations below are obtained:
\begin{equation} \nonumber
\begin{cases}
\frac{\sin{2\pi b}}{j+b}=2\pi\sum_{i=0}^{\infty}\frac{(-1)^i(2\pi(j+b))^{2i}}{(2i+1)!}\\

\frac{\sin{2\pi b}}{(j+b)^2}=\frac{2\pi}{j+b}-(2\pi)^2\sum_{i=0}^{\infty}\frac{(-1)^i(2\pi(j+b))^{2i+1}}{(2i+3)!}\\

\frac{\sin{2\pi b}}{(j+b)^3}=\frac{2\pi}{(j+b)^2}-(2\pi)^3\sum_{i=0}^{\infty}\frac{(-1)^i(2\pi(j+b))^{2i}}{(2i+3)!}\\

\frac{\sin{2\pi b}}{(j+b)^4}=\frac{2\pi}{(j+b)^3}-\frac{(2\pi)^3}{3!(j+b)}+(2\pi)^4\sum_{i=0}^{\infty}\frac{(-1)^i(2\pi(j+b))^{2i+1}}{(2i+5)!}\\

\frac{\sin{2\pi b}}{(j+b)^5}=\frac{2\pi}{(j+b)^4}-\frac{(2\pi)^3}{3!(j+b)^2}+(2\pi)^5\sum_{i=0}^{\infty}\frac{(-1)^i(2\pi(j+b))^{2i}}{(2i+5)!}\\

\frac{\sin{2\pi b}}{(j+b)^6}=\frac{2\pi}{(j+b)^5}-\frac{(2\pi)^3}{3!(j+b)^3}+\frac{(2\pi)^5}{5!(j+b)}-(2\pi)^6\sum_{i=0}^{\infty}\frac{(-1)^i(2\pi(j+b))^{2i+1}}{(2i+7)!}\\

\vdots 
\end{cases}
\end{equation}
\indent The formulae derived with these approaches break down for integer and half-integer $b$.

\subsection{Cosine approach}
Using the two systems of recurrence equations, we have more than one way to derive formulae for $\HP_k(n)$. Unless they are combined, one has a very complicated recurrence to determine the formula for $\HP_k(n)$. Let us go over the simplest way.\\

Setting $a=1$ without loss of generality, and intertwining the two systems of equations, the cosine-based $\HP_k(n)$ formula is,
\begin{equation} \nonumber
\sum_{j=1}^{n}\frac{1}{(j+b)^k}=-\frac{1}{2b^k}+\frac{1}{2(n+b)^k}-\frac{(2\pi)^k}{2}\int_{0}^{1}p_k(u)\left(\cos{2\pi(n+b)u}-\cos{2\pi bu}\right)\cot{\pi u}\,du \text{,}
\end{equation}
\noindent and the simplified recurrence equation on $p_k(u)$ is below:
\begin{equation} \nonumber
p_{k}(u)=\frac{(-1)^{\floor{(k+1)/2}}}{2\sin{\pi b}}\left(\frac{(1-u)^{k-1}}{(k-1)!}+\sum_{j=1}^{\floor{k/2}}\frac{(-1)^j p_{2j-1}(u)}{(k-2j+1)!}\right)\cdot
\begin{cases}
      -\frac{1}{\sin{\pi b}}, & \text{if }\ k \text{ is odd,}\\
      \frac{1}{\cos{\pi b}}, & \text{if }\ k \text{ is even}
\end{cases} 
\end{equation}
\indent To solve this equation, it is necessary to separate odd from even. For the odd, one can make $q_k(u)=p_{2k+1}(u)$, and hence for $k\ge 0$:
\begin{equation} \nonumber
q_{k}(u)=\frac{(-1)^{k}}{2(\sin{\pi b})^2}\left(\frac{(1-u)^{2k}}{(2k)!}-\sum_{j=0}^{k-1}\frac{(-1)^j q_j(u)}{(2k-2j)!}\right)
\end{equation}
\indent The recurrence for even is very similar and is a function of $q_k(u)$. The generating functions of odd and even $p_k(u)$ are therefore,
\begin{equation} \nonumber
f(x)=\frac{x\cos{x(1-u)}}{\cos{x}-\cos{2\pi b}} \text{ and } g(x)=\sin{x}\frac{x\cos{x(1-u)}}{\cos{x}-\cos{2\pi b}} \text{,}
\end{equation}
\noindent and the final formulae can be expressed as,
\begin{multline} \label{eq:cos_odd}
\small
\sum_{j=1}^{n}\frac{1}{(j+b)^{2k+1}}=-\frac{1}{2b^{2k+1}}+\frac{1}{2(n+b)^{2k+1}}\\-\frac{(2\pi)^{2k+1}}{2}\int_{0}^{1}\frac{f^{(2k+1)}(0)}{(2k+1)!}\left(\cos{2\pi(n+b)u}-\cos{2\pi b u}\right)\cot{\pi u}\,du \text{,}
\end{multline}
\noindent and
\begin{multline} \nonumber 
\small
\sum_{j=1}^{n}\frac{1}{(j+b)^{2k}}=-\frac{1}{2b^{2k}}+\frac{1}{2(n+b)^{2k}}\\-\frac{(2\pi)^{2k}}{2\sin{2\pi b}}\int_{0}^{1}\left(\frac{(-1)^k(1-u)^{2k-1}}{(2k-1)!}+\frac{g^{(2k)}(0)}{(2k)!}\right)\left(\cos{2\pi(n+b)u}-\cos{2\pi b u}\right)\cot{\pi u}\,du 
\end{multline}

\subsection{Sine approach} \label{sine}
Since the process for sine is analogous, only the final formulae are shown. Let $f(x)$ and $g(x)$ be the functions,
\begin{equation} \nonumber
f(x)=\frac{x\sin{x(1-u)}}{\cos{x}-\cos{2\pi b}} \text{ and } g(x)=\sin{x}\frac{x\sin{x(1-u)}}{\cos{x}-\cos{2\pi b}}
\end{equation}\\
\indent The even and odd $\HP_k(n)$ formulae can be expressed as a function of their derivatives by,
\begin{multline} \label{eq:sin_even}
\small
\sum_{j=1}^{n}\frac{1}{(j+b)^{2k}}=-\frac{1}{2b^{2k}}+\frac{1}{2(n+b)^{2k}}\\+\frac{(2\pi)^{2k}}{2}\int_{0}^{1}\frac{f^{(2k)}(0)}{(2k)!}\left(\sin{2\pi(n+b)u}-\sin{2\pi b u}\right)\cot{\pi u}\,du \text{,}
\end{multline}
\noindent and,
\begin{multline} \nonumber 
\small
\sum_{j=1}^{n}\frac{1}{(j+b)^{2k+1}}=-\frac{1}{2b^{2k+1}}+\frac{1}{2(n+b)^{2k+1}}\\+\frac{(2\pi)^{2k+1}}{2\sin{2\pi b}}\int_{0}^{1}\left(\frac{(-1)^k(1-u)^{2k}}{(2k)!}+\frac{g^{(2k+1)}(0)}{(2k+1)!}\right)\left(\sin{2\pi(n+b)u}-\sin{2\pi b u}\right)\cot{\pi u}\,du 
\end{multline}

\end{document}